\documentstyle[fullpage,maple,11pt]{article}
\title{Formal Power Series}%
\author{
	     Dominik Gruntz\\
       Institute for Scientific Computing\\ 
       ETH Z\"urich\\
       CH-8092 Z\"urich \\
       {\tt gruntz@inf.ethz.ch}
\and
       Wolfram Koepf\\
       Konrad-Zuse-Zentrum f\"ur Informationstechnik\\ 
       Heilbronner Str.\ 10\\ 
       D-10711 Berlin\\
       {\tt koepf@zib-berlin.de}
}
\date{December 28, 1993}
\begin{document}
\maketitle
\vspace*{-1cm}
\begin{center}
Konrad-Zuse-Zentrum Berlin (ZIB), Preprint SC 93-31, 1993
\end{center}
\newcommand{\heading}[1]{\section{#1}}

\newcommand{\Maple}{{\sc Maple}}
\newtheorem{theorem}{\sc Theorem}
\newtheorem{lemma}[theorem]{\sc Lemma}
\newtheorem{definition}{\sc Definition}
\newcommand{\notexists}{\raisebox{1pt}{$/$}\!\!\!\exists\,\,}

\newcommand{\funkdeff}[4]{\left\{\begin{array}{ccc}
                                 #1 && \mbox{\rm{if} $#2$ } \\[4mm]
                                 #3 && \mbox{\rm{if} $#4$ } 
                                 \end{array}
                          \right.}
\def\erf{\mathop{\rm erf}\nolimits}
\def\arccosh{\mathop{\rm arccosh}\nolimits}
\newenvironment{Program}%
  {\par \bgroup
   \parindent = 0pt
   \advance \linewidth by -2\leftmargin
   \parshape=1 \leftmargin \linewidth
   \def\>{\advance \leftmargin by 2em%
          \advance \linewidth by -2em%
          \parshape=1 \leftmargin \linewidth \ignorespaces}%
   \def\Cr{\par \endgroup \begingroup \leavevmode \ignorespaces}%
   \begingroup \leavevmode \ignorespaces}%
  {\par \endgroup
   \egroup}



\makeatletter
\newcounter{algorithm}
\def\thealgorithm{\@arabic\c@algorithm}
\def\fps@algorithm{\fps@figure}
\def\ftype@algorithm{4}
\def\ext@algorithm{loa}
\def\fnum@algorithm{Algorithm~\thealgorithm}
\def\algorithm{\@float{algorithm}}
\let\endalgorithm\end@float
\@namedef{algorithm*}{\@dblfloat{algorithm}}
\@namedef{endalgorithm*}{\end@dblfloat}
\makeatother


%





\newcommand{\C}{{\rm {\mbox{C{\llap{{\vrule height1.52ex}\kern.4em}}}}}}
\newcommand{\N} {{\rm {\mbox{\protect\makebox[.15em][l]{I}N}}}}
\newcommand{\Q} {{\rm {\mbox{Q{\llap{{\vrule height1.5ex}\kern.5em}}}}}}
\newcommand{\R} {{\rm {\mbox{\protect\makebox[.15em][l]{I}R}}}}
\newcommand{\D} {{\rm {\mbox{\protect\makebox[.15em][l]{I}D}}}}
\newcommand{\Z} {{\rm {\mbox{\protect\makebox[.2em][l]{\sf Z}\sf Z}}}}

\heading{Introduction}
\noindent
Formal Laurent-Puiseux series of the form
	\begin{equation}
	f(x)=\sum\limits_{k=k_0}^{\infty}a_{k}x^{k/n}
	\label{eq:formalLPS}
	\end{equation}
with coefficients $a_{k}\in\C\;(k\in\Z)$
are important in many branches of mathematics.
\Maple\ supports the computation of {\em truncated\/} series with its 
{\tt series} command, and through the {\tt powerseries} package
\cite{GruntzShare} infinite series are available. In the latter
case, the series is represented as a table of coefficients 
that have already been determined together with a function for computing
additional coefficients. This is known as {\em lazy evaluation\/}.
But these tools fail, if one is interested in an explicit formula
for the coefficients $a_k$.

In this article we will describe the \Maple\ implementation of
an algorithm presented in~\cite{Koe92}--\cite{Koeortho} which 
computes an {\em exact\/} formal power series (FPS) of a
given function. This procedure will 
enable the user to reproduce most of the results of the
extensive bibliography on series~\cite{Han}. We will give
an overview of the algorithm and then present some parts of it
in more detail. 

This package is available through the \Maple-share library
with the name {\tt FPS}.
We will flavor this procedure with the following example.
\begin{maple}
> FormalPowerSeries(sin(x), x=0);

                           infinity
                            -----       k  (2 k + 1)
                             \      (-1)  x
                              )     ----------------
                             /         (2 k + 1)!
                            -----
                            k = 0
\end{maple}

\heading{Preliminary Results}
\noindent
To deal with many special functions, it is a good idea to consider
the {\em (generalized) hypergeometric series\/} 
\begin{equation}\label{eq:hypergeometric_series}
_{p}F_{q}\left.\left(\begin{array}{cccc}
a_{1}&a_{2}&\cdots&a_{p}\\
b_{1}&b_{2}&\cdots&b_{q}\\
            \end{array}\right| x\right)
:=
\sum\limits_{k=0}^{\infty} \frac
{(a_{1})_{k}\cdot(a_{2})_{k}\cdots(a_{p})_{k}}
{(b_{1})_{k}\cdot(b_{2})_{k}\cdots(b_{q})_{k}\,k!}x^{k}
\label{eq:coefficientformula}
\end{equation}
where $(a)_{k}$ denotes the {\em Pochhammer symbol\/} 
(or {\em shifted factorial}) defined by
\[
(a)_{k}:=\funkdeff{1}{k=0}{a\cdot (a+1)\cdots (a+k-1)}{k\in\N}\;.
\]
We have defined this function in our package under the name
{\tt pochhammer}.

The coefficients $A_k$ of the hypergeometric series 
$\sum\limits_{k=0}^\infty A_k x^k$ are the unique solution of the special 
recurrence equation (RE)
\[
A_{k+1}:=\frac{(k+a_1)\cdot(k+a_2)\cdots (k+a_p)}{(k+b_1)\cdot(k+b_2)\cdots
(k+b_q)(k+1)}\cdot A_k\quad(k\in\N)
\]
with the initial condition
	\[
	A_0:=1\;.
	\]
Note that $\frac{A_{k+1}}{A_k}$ is rational in $k$.
Moreover if $\frac{A_{k+1}}{A_k}$ is a rational function $R(k)$
in the variable $k$ then the corresponding function $f$ is connected with a
hypergeometric series; i.e., if $k=-1$ is a pole of $R$, then
$f$ corresponds to
a hypergeometric series evaluated at some point $Ax$ (where $A$ is
the quotient of the leading coefficients of the numerator and the
denominator of $R$); whereas, if $k=-1$ is no pole of $R$, then $f$ may
be furthermore shifted by some factor $x^s\; (s\in\Z)$.

We further mention that the function $f$ corresponding to the 
hypergeometric series 
	\[
	f(x):=\;
	_{p}F_{q}\left.\left(\begin{array}{cccc}
	a_{1}&a_{2}&\cdots&a_{p}\\
	b_{1}&b_{2}&\cdots&b_{q}\\
             \end{array}\right| x\right)
	\]
satisfies the differential equation (DE)
\begin{equation}
\theta (\theta+b_1-1)\cdots (\theta+b_q-1)f=x(\theta+a_1)\cdots(\theta+a_p)f
\label{eq:hyperDE}
\end{equation}
where $\theta$ is the differential operator $x\frac{d}{dx}$ . An inspection
of the hypergeometric DE~(\ref{eq:hyperDE})
shows that it is of the form
\begin{equation}
\sum\limits_{j=0}^{Q} \sum\limits_{l=0}^{Q} c_{lj} x^l f^{(j)}=0
\label{eq:form}
\end{equation}
with certain constants $c_{lj}\in\C$ and $Q=\max\{p,q\}+1$.
Because of their importance in our development,
we call a DE of the form~(\ref{eq:form}), i.e.\ a homogeneous
linear DE with polynomial coefficients, {\em simple}.

We extend the considerations to 
{\em formal Laurent-Puiseux series} (LPS) with a representation
\begin{equation}
f:=\sum\limits_{k=k_{0}}^{\infty}a_{k}x^{k/n}
\quad\quad(a_{k_0}\neq 0)
\label{eq:FF}
\end{equation}
for some $k_0\in\Z$, and $n\in\N$. LPS are formal Laurent series,
evaluated at 
$\sqrt[n] x$. A formal Laurent series ($n=1$) is a shifted FPS, 
and corresponds to a meromorphic 
$f$ with a pole of order $-k_0$ at the origin. The number $n$ in
development~(\ref{eq:FF}) is called the {\em Puiseux number\/} of 
(the given representation of) $f$.

\begin{definition}
{\bf (Functions of hypergeometric type).}\ \
An LPS $f$ with representation~(\ref{eq:FF})
is called to be {\em of hypergeometric type\/} if its 
coefficients $a_k$ satisfy a RE of the form
\begin{eqnarray}
a_{k+m}&=&R(k)\,a_k\quad\:\mbox{for $k\geq k_0$}
\label{eq:HypergeomType}
\\
a_{k}&=&A_{k}\quad\quad\quad\;\mbox{for $k=k_0,k_0+1,\ldots,k_0+m-1$}
\nonumber
\end{eqnarray}
for some $m\in\N$, $A_{k}\in\C\;(k=k_0+1,k_0+2,\ldots,k_0+m-1)$,
$A_{k_0}\in\C\setminus\{0\}$, and some rational function $R$. The number $m$
is then called the {\em symmetry number} of (the given representation of) $f$.
A RE of type~(\ref{eq:HypergeomType}) is also called to be of hypergeometric
type.
\end{definition}

We want to emphasize that the above terminology of functions of
hypergeometric type is pretty more general than the terminology of
a generalized hypergeometric function. It covers e.g.\ the function 
$\sin x$ which
is {\em not} a generalized hypergeometric function as obviously no RE of 
type~(\ref{eq:HypergeomType}) with $m=1$ holds for its series coefficients. 
So $\sin{x}$ is not of hypergeometric type with symmetry number $1$; it is,
however, of hypergeometric type with symmetry number $2$.
A more difficult example of the same kind is the function $e^{\arcsin x}$
which is neither even nor odd, and nevertheless 
turns out to be of hypergeometric type with symmetry number $2$, too.
Also functions like $x^{-5}\,\sin{x}$ are covered by the given
approach.
Moreover the terminology covers composite functions like $\sin\sqrt
x$, which do not have a Laurent, but a Puiseux series development.
Each LPS with symmetry number $m$, and Puiseux number $n$,
can be represented as the sum of $nm$ 
shifted $m$-fold symmetric functions.

We remark further that one can
extend the definition of functions of hypergeometric type to
include also the functions of the form $\int{f}$ where $f$ is an
arbitrary LPS (see~\cite{Koe92}, Section 8).
Note that because of the logarithmic terms
these functions, in general, do not represent LPS.
  
In~\cite{Koe92} it has been proven, that each LPS of hypergeometric type satisfies
a simple DE, which is essential for our development.
Now assume, a function $f$ representing an LPS is given. In order to find the 
coefficient formula, it is a reasonable approach to search for its DE,
to transfer this DE into its equivalent RE, and you are done by an
adaption of the coefficient formula for the hypergeometric function
corresponding to transformations on $f$ which preserve its hypergeometric type.
Below, we will discuss these steps in detail.

The outline of the algorithm to produce a formal Laurent-Puiseux
series expansion of the function $f$ with respect to the variable 
$x$ around $x=0$ is shown in Algorithm~\ref{alg:FPS}. 
Note that series expansions around other points can easily
be reduced to this case.

\heading{Search of the DE}
\noindent
In this section we present the algorithm that searches for a simple DE
of degree $k$ for a given function $f$. We set up the
equation
\[
   f^{(k)}(x) + \sum_{j=0}^{k-1} A_j\,f^{(j)}(x) = 0
\]
and expand it. Then we collect the coefficients of all the
rationally dependent terms and equate them to zero. For testing
whether two terms are rationally dependent, we divide one by the
other and test whether the quotient is a rational function in $x$
or not. This is an easy and fast approach. Of course, we could also
use the Risch normalization procedure~\cite{Risch,Bronstein}
to generate the rationally independent
terms, but first this normalization is rather expensive and
only works for elementary functions and second,
our simplified approach will not lead to 
wrong results, it may at most happen that we miss a simpler
solution, which, in practice, rarely happens, however.
This procedure by its own is available under the name {\tt SimpleDE}.

\pagebreak\noindent
\begin{algorithm}[h]
  \begin{Program}
  {\tt FormalPowerSeries(f, x=0)}\Cr
  \> {\tt for} $k := 1$ {\tt to} $k_{max}$ {\tt do}\Cr
  \> \> Search a simple DE of degree $k$ of the form\Cr
\begin{equation}\label{eq:DEform}
\mbox{\tt DE := }f^{(k)}(x) + \sum\limits_{j=0}^{k-1} A_j\,f^{(j)}(x) = 0
\end{equation}
  \> \> \> where the $A_j$ are rational functions in $x$.\Cr
  \> \> {\tt if} the search was successful, {\tt then}\Cr
  \> \> \> Convert the {\tt DE} into a recurrence equation of the form\Cr
\begin{equation}\label{eq:REform}
\mbox{\tt RE := }\sum\limits_{j=0}^{M} p_{j} a_{k+j}=0
\end{equation}
  \> \> \> \> for the coefficients $a_k$, where $p_{j}$ are polynomials in $k$ and $M \in \N$\Cr
  \> \> \> {\tt if} the {\tt RE} only contains one or two summands {\tt then}\Cr
  \> \> \> \> $f$ is of {\em hypergeometric type\/} and the {\tt RE} can be solved \Cr
  \> \> \> {\tt elif} the {\tt DE} has constant coefficients {\tt then}\Cr
  \> \> \> \> $f$ is of {\em exp-like type} and the {\tt RE} can also be solved.\Cr
  \> \> \> {\tt fi}\Cr
  \> \> {\tt fi};\Cr
  \> \> {\tt if} $f^{(k)}(x)$ is a rational function in $x$, {\tt then}\Cr
  \> \> \> use the {\em rational Algorithm\/} and integrate the result $k$ times.\Cr
  \> \> {\tt fi}\Cr
  \> {\tt od}
  \end{Program}
  \caption{Algorithm {\tt FormalPowerSeries}}
  \label{alg:FPS}
\end{algorithm}

The resulting differential equation only depends on the form of 
the derivatives $f^{(j)}(x)$ (which of course must be known to 
the system). As an example we look at the Airy wave function {\tt Ai},
whose derivative presently (in \Maple) is given as
\begin{maple}
> diff(Ai(x),x);

                               1/2                   3/2
                            x 3    BesselK(2/3, 2/3 x   )
                      - 1/3 -----------------------------
                                          Pi
\end{maple}
The simple DE then becomes
\begin{maple}
> SimpleDE(Ai(x),x);

               /   3      \                      /   2      \
               |  d       |      2 /  d      \   |  d       |
               |----- F(x)| x - x  |---- F(x)| - |----- F(x)| = 0
               |   3      |        \ dx      /   |   2      |
               \ dx       /                      \ dx       /
\end{maple}
which is not the simple DE of lowest degree valid for {\tt Ai}. This happens
since the second derivative is not expressed in terms of {\tt Ai}
(and {\tt diff(Ai(x),x)}) itself (or in other words since {\tt Ai}
itself is not expressed in terms of {\tt BesselK}).
We may introduce a new function {\tt newAi} by defining its derivatives

\begin{maple}
> `diff/newAi`      := (e,x) -> diff(e,x)*newAiPrime(e):
> `diff/newAiPrime` := (e,x) -> diff(e,x)*e*newAi(e):
> SimpleDE(newAi(x),x);

                           /   2      \
                           |  d       |
                           |----- F(x)| - x F(x) = 0
                           |   2      |
                           \ dx       /
\end{maple}
and we get the expected differential equation for the Airy wave function
{\tt Ai}.

It may happen, that for a given function $f(x)$ a DE of degree $k$ exists,
but which has neither constant coefficients (which we call the
{\em explike case}) nor is the corresponding RE
of hypergeometric type and hence no closed form for the LPS can be
computed. What we can do in this situation is to look for a DE of
higher degree which then will have free parameters as from the
existence of a DE of degree $k$ follows the existence of families
of DEs of higher degree. These parameters can be set freely and
we can try to choose them in such a way, that we
can use our tools to compute a formal LPS, i.e.\ we either need
a RE of hypergeometric type or a DE with constant coefficients.

For the first case we convert the DE into the corresponding RE of
the form~(\ref{eq:REform}) and try to set all the coefficients but two
to zero. For example, let $f(x)=x^{-1}\,e^x\,\sin{x}$. We find DEs of degree 2 and 
of degree 3, but none of the corresponding REs can be solved. The RE which
corresponds to the DE of degree 4 has the form
\[
   \sum_{j=0}^4 p_j(k)\,a_{k+j} = 0
\]
where
\begin{eqnarray*}
p_0(k) & = & 2 A_2+4 A_3+4\\
p_1(k) & = & - 4 A_2 - 10 A_3 - 16 - 2 A_3 k - 2 A_2 k\\
p_2(k) & = & 18 A_3 + 5 A_2 k + 6 A_2 + 48 + 8 k + 6 A_3 k + A_2 k^2\\
p_3(k) & = & (4 + k) (A_3 k^2  + 2 A_3 k - 24 - 3 A_3)\\
p_4(k) & = & (k + 5) (4 + k) (k^2 + k + 6)\\
\end{eqnarray*}
%
%
The solution of the equations $p_1(k)=0, p_2(k)=0, p_3(k)=0$ 
(forcing that only two terms of the sum remain)
with respect to $A_2$ and $A_3$ is
\[
A_{2} = -8 \frac{k^2+5k+12}{k^3+4 k^2+k-6},\quad
A_{3} =    \frac{24}{k^2+2k-3}.
\]
This results in the following RE of hypergeometric type 
(after multiplying by $\frac{(k+2)\,(k+3)}{k^2+k+6}$)
\[
\left(k+5\right)\left(k+4\right)\left(k+3\right)
\left(k+2\right)\,a_{k+1}
+4\,a_k = 0
\]
with symmetry number $m=4$. Of course, we try to keep the symmetry number,
that is the number of resulting sums,
as small as possible. The final result for this example is

\begin{maple}
> FormalPowerSeries(exp(x)*sin(x)/x, x);

/infinity                          \   /infinity                              \
| -----       k   (- k)    k  (4 k)|   | -----       k   (- k)    k  (1 + 4 k)|
|  \      (-1)  64      256  x     |   |  \      (-1)  64      256  x         |
|   )     -------------------------| + |   )     -----------------------------|
|  /              (1 + 4 k)!       |   |  /           (1 + 4 k)! (2 k + 1)    |
| -----                            |   | -----                                |
\ k = 0                            /   \ k = 0                                /

       /infinity                               \
       | -----        k   (- k)    k  (4 k + 2)|
       |  \       (-1)  64      256  x         |
     + |   )     ------------------------------|
       |  /      (4 k + 3) (1 + 4 k)! (2 k + 1)|
       | -----                                 |
       \ k = 0                                 /
\end{maple}
Note that one of the four sums (with the powers $x^{4k+3}$) vanishes 
as a result of a vanishing initial coefficient.

If this step fails, then we try to choose the parameters such that the DE 
gets constant coefficients. Let us look at a rather similar example, 
$f(x)=x\,e^{x}\,\sin(2x)$. We again find DEs of second and third order,
but their corresponding REs cannot be solved. The DE of degree 4 has two free parameters
$A_2$ and $A_3$ as expected
\begin{eqnarray*}
\lefteqn{x^2\,\frac{d^4}{dx^4}f(x) +
   A_{3}{x}^{2}\,\frac{d^3}{dx^3}f(x) +
   A_{2}{x}^{2}\,\frac{d^2}{dx^2}f(x) +}\\
\quad&&
   \left(
      (A_3-2\,A_2+12)\,x^2+(-6\,A_3-2\,A_2+4)\,x
   \right)\,\frac{d}{dx}f(x) +\\
\quad&&
   \left(
      (10\,A_3+5\,A_2-5)x^2+(14\,A_3+2\,A_2+28)x +(6\,A_3+2\,A_2-4)
   \right)f(x)=0
\:.
\end{eqnarray*}
The DE will have constant coefficients, if 
both $A_2$ and $A_3$ are constants and if they meet the following equations 
\[
\begin{array}{rcc}
    6 A_3 + 2 A_2 -  4 & = & 0 \\
   14 A_3 + 2 A_2 + 28 & = & 0
\end{array}
\;.
\]
The solution of this system of equations is $A_2 = 14$ and $A_3 = -4$. If we
insert these values in the differential equation and divide by $x^2$, then we get
the DE
\[
   \frac{d^4}{dx^4}f(x) - 4\frac{d^3}{dx^3}f(x) + 14\frac{d^2}{dx^2}f(x)
   -20\frac{d}{dx}f(x) + 25f(x) = 0
\]
which has constant coefficients leading to a constant coefficient RE for 
$b_k$ given by $f(x)=\sum\limits\frac{b_k}{k!}\,x^k$ that can be solved:
\begin{maple}
> FormalPowerSeries(x*exp(x)*sin(2*x),x);

infinity
 -----   /        k 1/2                            k 1/2                   \
  \      |      (5 )    cos(k arctan(2)) k       (5 )    sin(k arctan(2)) k|  k
   )     |- 2/5 -------------------------- + 1/5 --------------------------| x
  /      \                  k!                               k!            /
 -----
 k = 0
\end{maple}

\heading{Conversion to the Recurrence Equation}
\noindent
As it has been proven in~\cite{Koe92}, this transformation is done
by the substitution
\begin{equation}
x^l f^{(j)}(x) \mapsto (k\!+\!1\!-\!l)_j \cdot a_{k+j-l}
\label{eq:rule}
\end{equation}
into the DE.

We give here a small \Maple\ procedure which performs this transformation.
We assume, that the $j^{th}$ derivative of $f(x)$ is represented by
the expression {\tt f(j)}. You may compare the procedure {\tt ConvertDEtoRE}
with the rule based solution presented in~\cite{Koe92}.

\begin{maple}
> ConvertDEtoRE := proc(de, f, x, a, k) local X, F, l, j;
>    if type(de,`+`) then
>       map(ConvertDEtoRE, de, f, x, a, k)
>    else
>       X := select(has, j*de, x);
>       F := select(has, j*de, f);
>       l := degree(X, x);
>       j := op(1,F);
>       de/X/F * pochhammer(k+1-l,j) * a(k+j-l)
>    fi
> end:
\end{maple}
The following example converts the left hand side of the DE for $e^x$ into
the corresponding left hand side of the RE for the coefficients $a_k$ of the
FPS of $e^x$.
\begin{maple}
> ConvertDEtoRE(F(1)-F(0), F, x, a, k);

                            (k + 1) a(k + 1) - a(k)
\end{maple}

The search of a DE and its conversion to the RE is 
directly available through the command {\tt SimpleRE}.
We see that {\tt newAi$(x^2)$} is of hypergeometric type
with symmetry number $m=6$, whereas the second RE is
not of hypergeometric type.
\begin{maple}
> SimpleRE(newAi(x^2), x);

                   (k - 1) (k + 1) a(k + 1) - 4 a(k - 5) = 0

> SimpleRE(x/(1-x-x^2), x);

             (1 - k) a(k) + (k - 1) a(k - 1) + (k - 1) a(k - 2) = 0
\end{maple}
The latter example is the generating function of the Fibonacci
numbers, and we get the expected RE. Note, that the common
factor $(k-1)$ ensures, that the RE holds $\forall\ k\in\Z$.

\heading{Solving a Recurrence Equation of Hypergeometric Type}
\noindent
If the recurrence equation of a function $f(x)$ is of the form
\begin{equation}\label{eq:RE1} 
   Q(k)\, a_{k+m} = P(k)\, a_k
\end{equation}
($P,Q$ polynomials)
then $f$ is of hypergeometric type and the corresponding series 
representation has symmetry number $m$. The explicit formula
for the coefficients can be found by the hypergeometric coefficient
formula~(\ref{eq:hypergeometric_series})
and some initial conditions. Based on the analysis of the polynomials $P(k)$
and $Q(k)$, we will convert the RE into one corresponding
to a {\em Taylor series\/} by applying a sequence of transformations
stated in the following lemma which preserve the hypergeometric type.

\begin{lemma}\label{Lemma:trafos}
Let $f$ be a formal Laurent series (FLS) of hypergeometric type with
representation~(\ref{eq:FF}),
whose coefficients $a_k$ satisfy a RE of the form~(\ref{eq:HypergeomType}), then
the following functions are of hypergeometric type, too. Their coefficients
$b_k$ satisfy a RE whose relation to the RE of $f$ is also given.
\[
\begin{array}{llllll}
\mbox{(a)} &\quad& x^n f(x)   &\quad& b_{k+m} = R(k-n)\,b_k & n \in \Z\\
\mbox{(b)} && f(A x)     && b_{k+m} = A^m\,R(k)\,b_k & A \in \C\\
\mbox{(c)} && f(x^n)     && b_{k+n\,m} = R(k/n)\,b_k & n \in \N\\
\mbox{(d)} && f(x^{1/m}) && b_{k+1} = R(k\,m)\,b_k &\\
\mbox{(e)} && f'(x)      && b_{k+m} = \frac{k+m+1}{k+1}\,R(k+1)\,b_k & \\
\end{array}
\]
\end{lemma}
For a proof of this lemma we refer to~\cite{Koe92}, Lemma 2.1 and Theorem
8.1.

First of all, we inspect the roots of $P(k)$ and $Q(k)$ of the
RE~(\ref{eq:RE1}). If there
are any rational roots, then we know that $f$ corresponds (possibly) to a
Puiseux series. In this case we transform $f$ to a function of Laurent
type by an application of transformation $(c)$, where $n$ is the
least common multiple of the denominators of all rational roots
of $P(k)$ and $Q(k)$. The FLS we get when we solve the transformed
RE can be transformed back to the LPS of $f$ by substituting
$x$ by $x^{1/n}$.

Let us now assume that $f$ can be expanded in a FLS. We reduce
this problem to solving the RE of a function which has a
FPS expansion. For that we remove the finite pole of $f$ at the
origin by multiplying $f$ with a suitable power of $x$
(transformation $(a)$). From the information of the RE
we can determine which power we have to use. Let us assume
that $f$ may be expanded in a Laurent series, i.e.\ that 
$\exists\ k_0\ :\ \forall\ k\leq k_0\ :\ a_k=0$. From these
known coefficients we can derive further ones using the 
given RE in the form
	\begin{equation}\label{eq:RE}
	a_{k+m} = \frac{P(k)}{Q(k)}\, a_k.
	\end{equation}
If we know that $a_k=0$ then also
$a_{k+m}=0$ provided that $Q(k) \neq 0$. Let $k_{min}$
be the smallest integer root of $Q(k)$. Consequently
$a_k=0\ \forall\  k < k_{min}+m$. From this it follows, that
	\[
	g = x^{-(k_{min}+m)} \, f
	\]
may be expanded into a FPS, i.e.\ $g$ has no pole at the origin,
given the assumption that $f$ may be expanded in a Laurent series. 
This latter assumption can
be tested by computing the limit of $g$ as $x$ tends to $0$. This
limit must be finite. (Since we also allow logarithmic singularities,
we test in fact whether the limit of $x\cdot g$ is 0.) If this is not
the case, then the assumption that $f$ may be expanded into a FLS 
is wrong and $f$ must have an essential singularity, i.e.\ 
$\notexists k_0\ :\ \forall\ k\leq k_0\ :\ a_k=0$. 
Otherwise the FPS of $g$ exists and can again easily be transformed
back to the FLS of $f$.

What we finally have to show is how to solve a RE which corresponds
to a given function $f$ which has a FPS expansion.
The RE~(\ref{eq:RE}) is valid $\forall\ k\ :\ Q(k)\neq 0$, especially
$\forall\ k>k_{max}$ where $k_{max}$ is the largest root of $Q(k)$.
Consequently we must determine $a_k$ for $k=0,1,\ldots,k_{max}+m$ and
have to solve the hypergeometric RE for $x^{-k}\,f(x^{1/m})$
using~(\ref{eq:hypergeometric_series}) for $k>k_{max}$.
We investigate now how this condition can be weakened.

The coefficient $a_k$ is given by the limit
$\lim\limits_{x\to 0} f^{(k)}(x)/k!$. 
Let us assume, that this limit is finite. Then the hypergeometric RE
can be solved in the case that $\forall\ j\geq 0\ :\ Q(k+j\,m)\neq 0$,
otherwise simply $a_k\,x^k$ is added to the
result.
Moreover note, that
\[
   \left(P(k)=0 \vee a_k=0 \right) \wedge 
   \forall\ j\geq 0\ :\ Q(k+j\,m)\neq0
   \wedge a_k \;\mbox{is finite}
   \Longrightarrow
   \forall\ j>0\ :\ a(k+j\,m) = 0
\]
and in this case the RE does not have to be solved and it is enough 
to add $a_k\,x^k$ to the result. The indices $k+j\,m, j>0$
no longer need to be considered.

If $a_k$ is infinite, then we found a logarithmic singularity which
we can remove by working with \[g=(x^{-k}\,f(x))'\] and by integrating
and shifting the resulting power series $S$. The RE of $g$ can be
obtained from the RE of $f$ by applying transformations $(a)$ and $(e)$.
The constant term which
we lose by the differentiation can be determined by computing
the limit \[\lim_{x\to 0} f(x)/x^k-\int S(x)\,dx.\]
Note that $g$ in general is a function with a FLS expansion and we first
must remove the pole to get a function with a FPS expansion. This is the
reason why we have chosen a recursive implementation of the
RE solver. The procedure {\tt hypergeomRsolve}
accepts as parameters $f(x), P(k), Q(k)$ and the symmetry number $m$.
Every application of a transformation of Lemma~\ref{Lemma:trafos} is 
nothing else but a recursive call of the RE solver.
One may inspect which steps the algorithm performs by assigning 
the variable {\tt infolevel[FormalPowerSeries]}. As an example we trace the
computation for $f(x)=x^{-1}\,\sin{\sqrt{x}}$. A DE of degree 2
is found whose corresponding RE is of hypergeometric type. From
the root of $2k+3$ it follows that the
Puiseux number is 2 and so transformation $(c)$ with $n=2$ is applied.
The resulting function is of Laurent type.
The smallest integer root of $Q(k)$ of the transformed RE is $-4$ and the 
symmetry number is $m=2$, hence we multiply the function with $x^2$
and adjust the RE accordingly. We end up with the function $\sin{x}$ whose
FPS can be computed directly. This result is then transformed back to the
LPS of $f(x)$ according to the two transformations we applied.

\begin{maple}
> infolevel[FormalPowerSeries] := 4:
> FormalPowerSeries(sin(sqrt(x))/x,x);
FPS/FPS:   looking for DE of degree   1
FPS/FPS:   looking for DE of degree   2
FPS/FPS:   DE of degree   2   found.
FPS/FPS:   DE = 
                     2
                  4 x  F''(x) + 10 x F'(x) + (2 + x) F(x) = 0

FPS/hypergeomRE:   RE is of hypergeometric type.
FPS/hypergeomRE:   Symmetry number m = 1
FPS/hypergeomRE:   RE: 2 (k + 2) (2 k + 3) a(k + 1) = - a(k)
FPS/hypergeomRE:   RE modified by k = 1/2*k
FPS/hypergeomRE:   => f := sin(x)/x^2

FPS/hypergeomRE:   RE is of hypergeometric type.
FPS/hypergeomRE:   Symmetry number m = 2
FPS/hypergeomRE:   RE: (k + 4) (k + 3) a(k + 2) = - a(k)
FPS/hypergeomRE:   working with x^2*f
FPS/hypergeomRE:   => f := sin(x)

FPS/hypergeomRE:   RE is of hypergeometric type.
FPS/hypergeomRE:   Symmetry number m = 2
FPS/hypergeomRE:   RE:  (k + 2) (k + 1) a(k + 2) = - a(k)
FPS/hypergeomRE:   RE valid for all k >=    0
FPS/hypergeomRE:   a(0) =   0
FPS/hypergeomRE:   a(2*j) = 0   for all j>0.
FPS/hypergeomRE:   a(1) =   1

                           infinity
                            -----       k  (k - 1/2)
                             \      (-1)  x
                              )     ----------------
                             /         (2 k + 1)!
                            -----
                            k = 0

> infolevel[FormalPowerSeries] := 1:
\end{maple}

\heading{Rational Algorithm}
\noindent
If the given function (or any of its derivatives) is rational in $x$ we can apply
the rational algorithm as described in Section 4 of~\cite{Koe92}.
First the complex partial fraction decomposition of $f(x)$
has to be calculated. Each term of the form $\frac{c}{(x-\alpha)^j}$
can be expanded by the binomial series whose coefficients are
\[a_k = \frac{(-1)^j\,c}{\alpha^{j+k}}
          \left(\begin{array}{c}j+k-1\\k\end{array}\right).
\]

\begin{maple}
> FormalPowerSeries(1/((x-1)^2*(x-2)),x);

                infinity
                 -----                  k                  k   k
                  \            (k! - 2 2  k! + 2 (1 + k)! 2 ) x
                   )     - 1/2 ---------------------------------
                  /                           k
                 -----                       2  k!
                 k = 0

> FormalPowerSeries((1+x+x^2+x^3)/((x-1)*(x-2)),x);

                             /infinity                   \
                             | -----           k        k|
                             |  \          (8 2  - 15) x |
                     x + 4 + |   )     1/2 --------------|
                             |  /                 k      |
                             | -----             2       |
                             \ k = 0                     /

> FormalPowerSeries(C/(B*A - A*x - B*x+x^2),x);

                          infinity
                           -----         k      k   k
                            \      C (A A  - B B ) x
                             )     ------------------
                            /                  k    k
                           -----    (A - B) B B  A A
                           k = 0
\end{maple}

To get the complex partial fraction decomposition we must factor
the denominator which may be rather complicated, hence
the rational algorithm to compute the full partial fraction expansion
presented in~\cite{Salvy} may be used. The following example uses this code. 
This method can be forced to be used, if the environment variable
\verb|_EnvExplicit| is set to false.

\begin{maple}
> FormalPowerSeries(1/(x^4+x+1),x);

      infinity
       -----   /  -----                                 2           3\
        \      |   \            27 + 64 alpha - 48 alpha  + 36 alpha |  k
         )     |    )     1/229 -------------------------------------| x
        /      |   /                              (1 + k)            |
       -----   |  -----                      alpha                   |
       k = 0   \alpha = 

                                 4           
\end{maple}

A closed form of the Fibonacci numbers can be derived by computing the
FPS of their generating function $x/(1-x-x^2)$.
\begin{maple}
> expand(FormalPowerSeries(x/(1-x-x^2),x));

                  /infinity                                            \
                  | -----   /           k                    k        \|
              1/2 |  \      |          x                    x         ||
       - 1/5 5    |   )     |- ----------------- + -------------------||
                  |  /      |        1/2       k                 1/2 k||
                  | -----   \  (1/2 5    - 1/2)    (- 1/2 - 1/2 5   ) /|
                  \ k = 0                                              /
\end{maple}
The {\tt expand} command converts the coefficient in the usual 
notation. If the factorization of the numerator is avoided, 
then the following result is obtained:
\begin{maple}
> _EnvExplicit := false:
> FormalPowerSeries(x/(1-x-x^2),x);

        infinity
         -----   /            -----                                \
          \      |             \                         alpha - 2 |  k
           )     |              )                - 1/5 ------------| x
          /      |             /                            (1 + k)|
         -----   |            -----                    alpha       |
         k = 0   |                            2                    |
                 \alpha = RootOf(- 1 + _Z + _Z )                   /

> _EnvExplicit := '_EnvExplicit':
\end{maple}

\heading{Special Functions, in Particular Orthogonal Polynomials}
\noindent
The algorithm has been extended to handle many special functions, in
particular orthogonal polynomials~\cite{Koeortho}.
Our implementation covers this approach.

We have seen, that the only precondition which must be met by a function to be covered by 
the algorithm is that its derivative is defined in terms of functions which also may be handled
by our algorithm. 

In the case of families of orthogonal polynomials we have the following 
special situation: The general derivative of an orthogonal polynomial
of degree $n$ can be defined in terms of the orthogonal polynomials
of degree $n$, and $n-1$. But on the other hand, furthermore a recurrence 
equation is known which allows to express the polynomial of degree $n$ in 
terms of the polynomials of degree $n-1$, and $n-2$. Combining these two 
facts, it is possible to find a second order DE for the general polynomial 
of degree $n$.
If the resulting RE is of hypergeometric type (which depends on the
point of expansion) the algorithm further yields an LPS representation
for the general polynomial of degree $n$.

Similarly if a function family $F(n,x)$ possesses a differentiation rule of the
form 
\[
\frac{\partial F(n,x)}{\partial x} =
p_0(n,x)\,F(n,x)+p_1(n,x)\,F(n-1,x)
\]
with rational expressions 
$p_0$, and $p_1$ (or a similar rule with $m$ rather than 2 
expressions on the right) then by the product and chain rules of 
differentiation the second derivative has the form
\[
\frac{\partial^2 F(n,x)}{\partial x^2} =
q_0(n,x)\,F(n,x)+q_1(n,x)\,F(n-1,x)+q_2(n,x)\,F(n-2,x)
\]
with rational expressions $q_0, q_1$, and $q_2$, and all higher 
derivatives of $F(n,x)$ obtain similar representations. Each differentiation
increases the number of representing expressions by one. However, if we 
further know a recurrence equation of the form
\[
F(n,x)=r_1(n,x)\,F(n-1,x)+r_2(n,x)\,F(n-2,x)
\]
with rational expressions $r_1$, and $r_2$ (or a similar equation with
$m$ rather than 2 expressions on the right) then a recursive application
of this equation can be used to simplify each combination of derivatives 
of $F(n,x)$ to a sum of $2$ (or $m$, respectively) rationally independent ones.

To give an example, we consider the 
Fibonacci polynomials. The recurrence equation of the 
family of Fibonacci polynomials $F_n(x)$ is
\begin{eqnarray*}
F_n(x) &=& x\,F_{n-1}(x) + F_{n-2}(x)\\
F_0(x) &=& 0\\
F_1(x) &=& 1.
\end{eqnarray*}
We can teach our procedure to use this recurrence equation by
assigning the table {\tt FPSRecursion}. The second index specifies
the number of arguments of the function family.
\begin{maple}
> FPSRecursion[Fibonacci, 2] := (n,x) -> x*Fibonacci(n-1,x) + Fibonacci(n-2,x):
\end{maple}
The derivative rule is given by
\[
\frac{\partial F_n(x)}{\partial x} =
   \frac{(n-1)\,x\,F_n(x)+2\,n\,F_{n-1}(x)}{x^2+4}.
\]
and written in Maple
\begin{maple}
> `diff/Fibonacci` := proc(n, e, x)
>     diff(e,x) *((n-1)*e*Fibonacci(n,e)+2*n*Fibonacci(n-1,e))/(e^2+4)
> end:
\end{maple}
This rule has been derived as follows. First, an explicit
formula of $F_n(x)$ has been computed using our algorithm
to generate the FPS of the generating function
$\sum\limits_{n=0}^\infty F_n(x)\,t^n=\frac{t}{1-x\,t-t^2}$
of the Fibonacci polynomials that is an easy consequence of the recurrence 
equation:
\begin{maple}
> FormalPowerSeries(t/(1-x*t-t^2), t);

                             2     1/2 k                    2     1/2 k   k
         (- (- 1/2 x - 1/2 (x  + 4)   )  + (- 1/2 x + 1/2 (x  + 4)   ) ) t
 Sum(- ---------------------------------------------------------------------,
                        2     1/2 k                  2     1/2 k   2     1/2
       (- 1/2 x - 1/2 (x  + 4)   )  (- 1/2 x + 1/2 (x  + 4)   )  (x  + 4)

     k = 0 .. infinity)
\end{maple}
After some simplifications, the coefficient of this FPS, i.e.\ $F_n(x)$ 
has the following form:
\begin{maple}
> F := ((1/2*x+1/2*(x^2+4)^(1/2))^n-(1/2*x-1/2*(x^2+4)^(1/2))^n)/(x^2+4)^(1/2);

                              2     1/2 n                  2     1/2 n
               (1/2 x + 1/2 (x  + 4)   )  - (1/2 x - 1/2 (x  + 4)   )
          F := -------------------------------------------------------
                                       2     1/2
                                     (x  + 4)
\end{maple}
We now make the following ansatz for the derivative rule, namely
\[F_n(x)' = a\,F_n(x) + b\,F_{n-1}(x)\] and try to 
solve this equation for the unknown parameters $a$, and $b$:
\begin{maple}
> eq := diff(F,x) - (a*F + b*subs(n=n-1,F)):
> eq := numer(normal(eq, expanded)):
> indets(eq);
                            2     1/2                 2     1/2 n
             {n, a, b, x, (x  + 4)   , (1/2 x + 1/2 (x  + 4)   ) ,

                                2     1/2 n    2     3/2
                 (1/2 x - 1/2 (x  + 4)   ) , (x  + 4)   }

> solve({coeffs(eq, {"[5..8]})}, {a,b});

                              x (n - 1)           n
                         {a = ---------, b = 2 ------}
                                 2              2
                                x  + 4         x  + 4

> assign(");
> a*Fibonacci(n,x)+b*Fibonacci(n-1,x);

              x (n - 1) Fibonacci(n, x)     n Fibonacci(n - 1, x)
              ------------------------- + 2 ---------------------
                         2                           2
                        x  + 4                      x  + 4
\end{maple}
Note that, again, the above procedure essentially equates the coefficients
of the rationally independent terms of our setting to zero.

The algorithm trying to find a DE computes the first and the second derivative
of $F_n(x)$, expresses all occurrences of $F_n(x)$ in terms of $F_{n-1}(x)$ and 
$F_{n-2}(x)$, and finally returns the following solution:
\begin{maple}
> SimpleDE(Fibonacci(n,x),x);

                /   2      \
         2      |  d       |                              /  d      \
       (x  + 4) |----- F(x)| - (n - 1) (n + 1) F(x) + 3 x |---- F(x)| = 0
                |   2      |                              \ dx      /
                \ dx       /
\end{maple}
If we declare the initial value for $x=0$ which is 0 for even $n$ and 1 for odd
$n$ (which follows from the recurrence equation for $x=0$: $F_n(0)=F_{n-2}(0)$
and the initial conditions for $n=0$ and $n=1$) we can compute the formal
power series of the Fibonacci polynomial.
\begin{maple}
> Fibonacci := proc(n,x)
>    if x=0 then sin(n*Pi/2)^2
>    elif n=0 then 0
>    elif n=1 then 1
>    else 'procname(args)'
>    fi
> end:
> Fib := FormalPowerSeries(Fibonacci(n,x), x);

                    2
Fib := sin(1/2 n Pi)  (n - 1) (n + 1)

/infinity                                                     \
| -----       k                                          (2 k)|
|  \      (-1)  (- 1/2 n + 1/2 + k)! (1/2 n + 1/2 + k)! x     |
|   )     ----------------------------------------------------|/
|  /               (n - 2 k - 1) (n + 2 k + 1) (2 k)!         |
| -----                                                       |
\ k = 0                                                       /

((- 1/2 n + 1/2)! (1/2 n + 1/2)!)

                        /infinity                                             \
                        | -----       k                              (2 k + 1)|
                      2 |  \      (-1)  (- 1/2 n + k)! (1/2 n + k)! x         |
       n cos(1/2 n Pi)  |   )     --------------------------------------------|
                        |  /                       (2 k + 1)!                 |
                        | -----                                               |
                        \ k = 0                                               /
 + 1/2 ------------------------------------------------------------------------
                                  (- 1/2 n)! (1/2 n)!
\end{maple}
Note, that some of the factorials may have negative arguments,
and hence the limits of the coefficients must be considered.


Let's for example compute $F_{10}(x)$ and $F_{11}(x)$
which are polynomials of degree 9 and 10 respectively.
For $n=10=$ even we must only consider the second term of the
solution {\tt Fib} for $k$ up to $n/2-1=4$ (as we shall show soon).
Similarly for $n=11=$ odd we only consider the even coefficients
of {\tt Fib} for $k$ up to $(n-1)/2=5$.

\begin{maple}
> limit(subs(Sum=sum, infinity=4, op(2,Fib)), n=10);

                         9      7       5       3
                        x  + 8 x  + 21 x  + 20 x  + 5 x

> limit(subs(Sum=sum, infinity=5, op(1,Fib)), n=11);

                      10      8       6       4       2
                     x   + 9 x  + 28 x  + 35 x  + 15 x  + 1
\end{maple}

In the second term of the above solution, the form of the odd
coefficients are
\begin{eqnarray*}
a_{2k+1} &= &
\frac{n\,\cos^2\left(\frac{n}{2}\pi\right)}{2}
\frac{(-1)^k}{(2k+1)!}
\frac{\left(\frac{n}{2}+k\right)!}
                     {\left(\frac{n}{2}\right)!}
\frac{\left(-\frac{n}{2}+k\right)!}
                     {\left(-\frac{n}{2}\right)!}
\\&=&
\frac{n\,\cos^2\left(\frac{n}{2}\pi\right)}{2}
\frac{(-1)^k}{(2k+1)!}
\frac{\Gamma\left(1+\frac{n}{2}+k\right)}
		{\Gamma\left(1+\frac{n}{2}\right)}
\frac{\Gamma\left(1-\left(\frac{n}{2}-k\right)\right)}
                {\Gamma\left(1-\frac{n}{2}\right)}
\\&=&
\frac{n\,\cos^2\left(\frac{n}{2}\pi\right)}{2}
\frac{(-1)^k}{(2k+1)!}
\frac{\Gamma\left(1+\frac{n}{2}+k\right)}
                {\Gamma\left(1+\frac{n}{2}\right)}
\frac{\Gamma\left(\frac{n}{2}\right)}
                {\Gamma\left(\frac{n}{2}-k\right)}
\frac{\sin\left(\frac{n}{2}\pi\right)}
		{\sin\left(\pi\left(\frac{n}{2}-k\right)\right)}
\\&=&
\frac{\cos^2\left(\frac{n}{2}\pi\right)\,(-1)^k}{(2k+1)!}
\frac{\Gamma\left(1+\frac{n}{2}+k\right)}
                {\Gamma\left(\frac{n}{2}-k\right)}
\frac{\sin\left(\frac{n}{2}\pi\right)}
                {\sin\left(\pi\left(\frac{n}{2}-k\right)\right)}
\\&=&
\frac{\cos^2\left(\frac{n}{2}\pi\right)}{(2k+1)!}
\frac{\Gamma\left(1+\frac{n}{2}+k\right)}
                {\Gamma\left(\frac{n}{2}-k\right)}
=
\funkdeff{\frac{1}{(2k+1)!}
\frac{\left(\frac{n}{2}+k\right)!}{\left(\frac{n}{2}-k-1\right)!}
}{n \mbox{ is even}}{0}{n \mbox{ is odd}}
\end{eqnarray*}
where we used the identities 
\[
\Gamma(x+1)=x\Gamma(x)=x!\;, 
\]
\[
\Gamma(x)\Gamma(1-x)=\frac{\pi}{\sin(\pi x)}
\;,
\]
and
\[
\sin(a+b)=\sin a\cos b+\cos a \sin b\;.
\]
Similarly we get for the even coefficients
\begin{eqnarray*}
a_{2k} &=&
\frac{\sin^2\left(\frac{n}{2}\pi\right)\,(-1)^k}{(2k)!}
\frac{\left(\frac{1}{2}(n+1)+k\right)!}{n+2k+1}
\frac{n+1}{\left(\frac{1}{2}(n+1)\right)!}
\frac{\left(-\frac{1}{2}n+\frac{1}{2}+k\right)!}{n-2k-1}
\frac{n-1}{\left(-\frac{1}{2}n+\frac{1}{2}\right)!}
\\&=&
\frac{\sin^2\left(\frac{n}{2}\pi\right)\,(-1)^k}{(2k)!}
\frac{\left(\frac{1}{2}(n-1)+k\right)!}{\left(\frac{1}{2}(n-1)\right)!}
\frac{\left(-\frac{1}{2}(n+1)+k\right)!}{\left(-\frac{1}{2}(n+1)\right)!}
\\&=&
\frac{\sin^2\left(\frac{n}{2}\pi\right)\,(-1)^k}{(2k)!}
\frac{\Gamma\left(1+\frac{1}{2}(n-1)+k\right)}
	{\Gamma\left(1+\frac{1}{2}(n-1)\right)}
\frac{\Gamma\left(1-\frac{1}{2}(n+1)+k\right)}
        {\Gamma\left(1-\frac{1}{2}(n+1)\right)}
\\&=&
\frac{\sin^2\left(\frac{n}{2}\pi\right)\,(-1)^k}{(2k)!}
\frac{\Gamma\left(1+\frac{1}{2}(n-1)+k\right)}
        {\Gamma\left(1+\frac{1}{2}(n-1)\right)}
\frac{\Gamma\left(\frac{1}{2}(n+1)\right)}
	{\Gamma\left(\frac{1}{2}(n+1)-k\right)}
\frac{\sin\left(\frac{\pi}{2}(n+1)\right)}
        {\sin\left(\left(\frac{1}{2}(n+1)-k\right)\pi\right)}
\\&=&
\frac{\sin^2\left(\frac{n}{2}\pi\right)}{(2k)!}
\frac{\Gamma\left(1+\frac{1}{2}(n-1)+k\right)}
        {\Gamma\left(\frac{1}{2}(n+1)-k\right)}
=
\funkdeff{\frac{1}{(2k)!}
\frac{\left(\frac{1}{2}(n-1)+k\right)!}{\left(\frac{1}{2}(n-1)-k\right)!}
}{n \mbox{ is odd}}{0}{n \mbox{ is even}}
\end{eqnarray*}
Note that $a_{2k+1}=0$ for $k\geq\frac{n}{2}$, and $a_{2k}=0$ 
for $k\geq\frac{n+1}{2}$.

If we make use of this additional information we end up with the
following closed formula for the general Fibonacci polynomial
$F_n(x)$. 
\[
F_n(x) =\funkdeff{
\sum\limits_{k=0}^{n/2-1}\frac{1}{(2k+1)!}
\frac{\left(\frac{n}{2}+k\right)!}{\left(\frac{n}{2}-k-1\right)!}x^{2k+1}
}{n \mbox{ is even}}{
\sum\limits_{k=0}^{(n-1)/2}\frac{1}{(2k)!}
\frac{\left(\frac{1}{2}(n-1)+k\right)!}{\left(\frac{1}{2}(n-1)-k\right)!}x^{2k}
}{n \mbox{ is odd}}
\;.
\]
%
%

Our implementation covers derivative rules and uses recurrence
equations for the following families of
special functions: the Fibonacci polynomials
{\tt Fibonacci(n,x)}, the Bessel functions {\tt BesselJ(n,x)},
{\tt BesselY(n,x)}, {\tt BesselI(n,x)}, and {\tt BesselK(n,x)}
(see~\cite{AS}, (9.1) and (9.6)),\linebreak
the Hankel functions {\tt Hankel1(n,x)}, and {\tt Hankel2(n,x)}
(see~\cite{AS}, (9.1)),
the Kummer\linebreak functions
{\tt KummerM(a,b,x)}, and {\tt KummerU(a,b,x)} (see~\cite{AS}, (13.4)),
the Whittaker functions \linebreak
{\tt WhittakerM(n,m,x)}, 
and {\tt WhittakerW(n,m,x)} (see~\cite{AS}, (13.4)),
the associated Legendre functions {\tt LegendreP(a,b,x)}, 
and {\tt LegendreQ(a,b,x)} (see~\cite{AS}, (8.5)),
the orthogonal polynomials \linebreak
{\tt JacobiP(n,a,b,x)}, 
{\tt GegenbauerC(n,a,x)}, {\tt ChebyshevT(n,x)}, {\tt ChebyshevU(n,x)},
{\tt LegendreP(n,x)}, {\tt LaguerreL(n,a,x)}, and {\tt HermiteH(n,x)}
(see~\cite{AS}, (22.7) and (22.8)), 
and the iterated integrals of the complementary error function
{\tt erfc(n,x)} (see~\cite{AS}, (7.2)).

As orthogonal polynomials are polynomials, for each fixed number $n$ there 
is a simple DE of order one. This is the reason why we use different names
from those in the packages {\tt orthopoly} or {\tt combinat}, as otherwise
for fixed $n$ evaluation occurs, and a first order DE is created which does 
not possess the structure of the polynomial system, and second, the
derivative of {\tt orthopoly[T]} can not be defined. For example
\begin{maple}
> SimpleDE(LaguerreL(3,a,x),x);

             /   2      \
             |  d       |                          /  d      \
             |----- F(x)| x + 3 F(x) + (a + 1 - x) |---- F(x)| = 0
             |   2      |                          \ dx      /
             \ dx       /
\end{maple}
generates the second order DE which structurally characterizes
the third Laguerre polynomial, 
whereas with the Laguerre polynomial out of the {\tt orthopoly} package
we get
\begin{maple}
> SimpleDE(orthopoly[L](3,a,x),x);

                          2              2
   (18 + 15 a - 18 x + 3 a  - 6 a x + 3 x ) F(x) +

                             2               2    3      2          2    3
       (6 + 11 a - 18 x + 6 a  - 15 a x + 9 x  + a  - 3 a  x + 3 a x  - x )

       /  d      \
       |---- F(x)| = 0
       \ dx      /
\end{maple}
We note that by a general result the sum, product, and composition with
rational functions and rational powers of functions satisfying the
type of DE considered, inherit this property~\cite{Koediffgl}, 
so that the algorithm is capable to generate the DE 
of a large class of functions. Here we give some more examples:

\begin{maple}
> SimpleDE(exp(x/2)*WhittakerW(a,b,1/x),x);

       4      3      2  2    2                        3         /  d      \
     (x  - 4 x  - 4 x  b  + x  + 4 a x - 1) F(x) - 4 x  (x - 2) |---- F(x)|
                                                                \ dx      /

              /   2      \
              |  d       |  4
          + 4 |----- F(x)| x  = 0
              |   2      |
              \ dx       /

> SimpleDE(LegendreP(n,1/(1-x)),x);

                    /   2      \
                  2 |  d       |                             3 /  d      \
 x (x - 2) (x - 1)  |----- F(x)| + n (n + 1) F(x) + 2 (x - 1)  |---- F(x)| = 0
                    |   2      |                               \ dx      /
                    \ dx       /

> FormalPowerSeries(JacobiP(n,a,b,x),x=1);
FPS/hypergeomRE:   provided that    -1    <=    min(-1,-a-1)

         binomial(a + n, n) n a!

             /infinity                                               \
             | -----                                   (- k)        k|
             |  \      (- n + k)! (b + a + n + k)! (-2)      (x - 1) |
             |   )     ----------------------------------------------|/
             |  /                    (n - k) (a + k)! k!             |
             | -----                                                 |
             \ k = 0                                                 /

             ((- n)! (b + a + n)!)
\end{maple}
Note, that here again for $(- n + k)!/(- n)!$ the corresponding limits must be
considered.

\heading{Asymptotic Series}
\noindent
The same algorithm can also be used to compute asymptotic expansions.
Note, that we only look for asymptotic series of the form of a
Laurent-Puiseux series. Such series are unique, a property
which is not given for general asymptotic expansions. As a consequence
the results we compute may differ from the truncated asymptotic
expansions returned by the \Maple\ command {\tt asympt}. 

One special thing of Laurent-Puiseux asymptotic expansions 
is, that they are only valid as long as the indeterminate 
approaches the expansion point from one side. If the expansion
point is not $\infty$, then one may specify with an option 
{\tt right} or {\tt left} from which side one approaches the
expansion point.

\begin{maple}
> FormalPowerSeries(erf(x), x=infinity);

                                       1

> FormalPowerSeries(arctan(1/x), x=0, right);

                               /infinity                 \
                               | -----       k  (2 k + 1)|
                               |  \      (-1)  x         |
                      1/2 Pi - |   )     ----------------|
                               |  /           2 k + 1    |
                               | -----                   |
                               \ k = 0                   /

> FormalPowerSeries(exp(x), x=-infinity);

                                       0

> FormalPowerSeries(exp(x), x=infinity):
FPS/FPS:   ERROR:    essential singularity

> FormalPowerSeries(exp(x)*Ei(-x) + exp(-x)*Ei(x), x=infinity);

                       /infinity                          \
                       | -----                            |
                       |  \                      (2 k + 2)|
                     2 |   )     (1 + 2 k)! (1/x)         |
                       |  /                               |
                       | -----                            |
                       \ k = 0                            /

> FormalPowerSeries(exp(x)*(1-erf(sqrt(x))), x=infinity);

                infinity
                 -----       (- k)         (- k)      (1/2 + k)
                  \      (-1)      (2 k)! 4      (1/x)
                   )     --------------------------------------
                  /                        k!
                 -----
                 k = 0
                -----------------------------------------------
                                       1/2
                                     Pi
\end{maple}
By a plot of {\tt arctan(1/x)}, e.g, one may realize that, indeed, the resulting
series representation is one sided, only.

\heading{Examples}
\noindent
In this section we present some more results generated with the procedure
{\tt FormalPowerSeries}.

%
%
First, we consider the FPS of the
Airy wave function {\tt Ai}. Since we use our own definition of the
derivatives, we also must define the initial values for $x=0$
(see e.g.~\cite{AS}, (10.4.4)--(10.4.5)).
\begin{maple}
> newAi(0)      :=  1/3^(2/3)/GAMMA(2/3):
> newAiPrime(0) := -1/3^(1/3)/GAMMA(1/3):
> FormalPowerSeries(newAi(x),x);

                      /infinity                      \
                      | -----        (- k)  (3 k)    |
                  1/3 |  \          9      x         |
                 3    |   )     ---------------------|
                      |  /      pochhammer(2/3, k) k!|
                      | -----                        |
                      \ k = 0                        /
             1/3 -------------------------------------
                               GAMMA(2/3)

                                     /infinity                      \
                                     | -----      (- k)  (3 k + 1)  |
                      1/6            |  \        9      x           |
                     3    GAMMA(2/3) |   )     ---------------------|
                                     |  /      pochhammer(4/3, k) k!|
                                     | -----                        |
                                     \ k = 0                        /
               - 1/2 ------------------------------------------------
                                            Pi
\end{maple}
Note, that \Maple\ is not yet able to compute even a truncated
power series of {\tt Ai(x)}. The same holds obviously for the
following example as long as the parameter $a$ is left as an unknown.
\begin{maple}
> FormalPowerSeries(x^a*sin(x^2),x);

                         infinity
                          -----       k  (4 k + 2 + a)
                           \      (-1)  x
                            )     --------------------
                           /           (1 + 2 k)!
                          -----
                          k = 0
\end{maple}

The next two examples are interesting results which may be unexpected.
It turns out, that both, $e^{\arccosh{x}}$ and $e^{\arccos{x}}$ are of
hypergeometric type.
\begin{maple}
> FormalPowerSeries(exp(arccos(x)),x);

                /         /  k - 1              \          \
                |         | --------'           |          |
                |         |'  |  |      2       |  k  (2 k)|
                |         |   |  |    (j  + 1/4)| 4  x     |
                |infinity |   |  |              |          |
                | -----   |   |  |              |          |
                |  \      \  j = 0              /          |
    exp(1/2 Pi) |   )     ---------------------------------|
                |  /                    (2 k)!             |
                | -----                                    |
                \ k = 0                                    /

                       /         /  k - 1                  \              \
                       |         | --------'               |              |
                       |         |'  |  |                2 |  k  (2 k + 1)|
                       |         |   |  |    (1/2 + j + j )| 4  x         |
                       |infinity |   |  |                  |              |
                       | -----   |   |  |                  |              |
                       |  \      \  j = 0                  /              |
         - exp(1/2 Pi) |   )     -----------------------------------------|
                       |  /                      (2 k + 1)!               |
                       | -----                                            |
                       \ k = 0                                            /

> FormalPowerSeries(exp(arccosh(x)),x);

                        /infinity                     \
                        | -----    (- k)         (2 k)|
                        |  \      4      (2 k)! x     |
                    - I |   )     --------------------| + x
                        |  /             2            |
                        | -----      (k!)  (2 k - 1)  |
                        \ k = 0                       /
\end{maple}
If we substitute {\tt Sum} (the inert form of summation)
by {\tt sum} which tries to get a closed formula, then we get the well known
identity $e^{\arccosh{x}}=x+\sqrt{x^2-1}$:
\begin{maple}
> eval(subs(Sum=sum,"));
                                       2 1/2
                               I (1 - x )    + x
\end{maple}
In fact, both $e^{\arccosh{x}}$ and $e^{\arccos{x}}$ are the special cases
$a=1$ and $a=\sqrt{-1}$ of the function $(x+\sqrt{x^2-1})^a$.

The following example has a logarithmic singularity at $x=0$ and hence
the FPS of $f'$ is computed.
\begin{maple}
> FormalPowerSeries(arcsech(x), x, real);

                              /infinity                             \
                              | -----               (- k)  (2 k + 2)|
                              |  \      (1 + 2 k)! 4      x         |
          ln(2) - ln(x) - 1/2 |   )     ----------------------------|
                              |  /             2                    |
                              | -----      (k!)  (k + 1) (2 k + 2)  |
                              \ k = 0                               /
\end{maple}

The next two examples are FPS expansions of two definite integrals, namely the
polylogarithm function and $\int\limits_0^x \erf(t)/t\,dt$. Note that for the
second integral we have used the inert function of {\tt Int} which
is only a placeholder and does not try to compute a closed form.
\begin{maple}
> FormalPowerSeries(dilog(1-x), x);

                               infinity
                                -----    (k + 1)
                                 \      x
                                  )     --------
                                 /             2
                                -----   (k + 1)
                                k = 0

> FormalPowerSeries(Int(erf(t)/t,t=0..x), x);

                            infinity
                             -----       k  (2 k + 1)
                              \      (-1)  x
                               )     ----------------
                              /                    2
                             -----     k! (2 k + 1)
                             k = 0
                          2 -------------------------
                                        1/2
                                      Pi
\end{maple}

The next two examples are FPS of functions which contain orthogonal
polynomials.
\begin{maple}
> FormalPowerSeries(exp(-x)*LaguerreL(n,a,x),x);
FPS/hypergeomRE:   provided that    -1    <=    min(-1,-a-1)

                                   /infinity                      \
                                   | -----       k               k|
                                   |  \      (-1)  (n + a + k)! x |
             binomial(n + a, n) a! |   )     ---------------------|
                                   |  /           (a + k)! k!     |
                                   | -----                        |
                                   \ k = 0                        /
             ------------------------------------------------------
                                    (n + a)!

> FormalPowerSeries(exp(-x^2)*HermiteH(n,x),x);

                            /infinity                                \
                            | -----                          k  (2 k)|
                            |  \      (1/2 n + 1/2 + k)! (-4)  x     |
     cos(1/2 n Pi) (n + 1)! |   )     -------------------------------|
                            |  /            (n + 1 + 2 k) (2 k)!     |
                            | -----                                  |
                            \ k = 0                                  /
     -----------------------------------------------------------------
                          (1/2 n)! (1/2 n + 1/2)!

                                   /infinity                              \
                                   | -----       k               (2 k + 1)|
                                   |  \      (-4)  (1/2 n + k)! x         |
            sin(1/2 n Pi) (n + 1)! |   )     -----------------------------|
                                   |  /                (2 k + 1)!         |
                                   | -----                                |
                                   \ k = 0                                /
          + ---------------------------------------------------------------
                                (1/2 n + 1/2)! (1/2 n)!
\end{maple}

As we have seen above, we have three tools available to compute a formal
power series representation of a given function, and it may happen for
some examples, that several of these tools may be applied. Normally the
solutions of a RE of hypergeometric type leads to the simplest results,
but this is not true in general. Hence we added the possibility for the
user to choose a method by adding an additional argument which is either
{\tt hypergeometric}, {\tt explike} or {\tt rational}.

\begin{maple}
> f := x*arctan(x)-1/2*ln(1+x^2):
> FormalPowerSeries(f, x, hypergeometric);

                            /infinity                  \
                            | -----        k  (2 k + 2)|
                            |  \       (-1)  x         |
                        1/2 |   )     -----------------|
                            |  /      (k + 1) (2 k + 1)|
                            | -----                    |
                            \ k = 0                    /

> FormalPowerSeries(f, x, rational);

                     infinity
                      -----               k    k   (k + 2)
                       \            ((- I)  + I ) x
                        )     1/2 -------------------------
                       /           k      k
                      -----       I  (- I)  (k + 1) (k + 2)
                      k = 0

> f :=sin(x)*exp(x):
> FormalPowerSeries(f, x, explike);

                       infinity
                        -----     k 1/2                k
                         \      (2 )    sin(1/4 k Pi) x
                          )     ------------------------
                         /                 k!
                        -----
                        k = 0

> FormalPowerSeries(f, x, hypergeometric);

                /infinity                              \
                | -----       k   (- k)    k  (4 k + 1)|
                |  \      (-1)  64      256  x         |
                |   )     -----------------------------|
                |  /                (4 k + 1)!         |
                | -----                                |
                \ k = 0                                /

                       /infinity                              \
                       | -----       k   (- k)    k  (4 k + 2)|
                       |  \      (-1)  64      256  x         |
                     + |   )     -----------------------------|
                       |  /           (2 k + 1) (4 k + 1)!    |
                       | -----                                |
                       \ k = 0                                /

                       /infinity                               \
                       | -----        k   (- k)    k  (4 k + 3)|
                       |  \       (-1)  64      256  x         |
                     + |   )     ------------------------------|
                       |  /      (2 k + 1) (4 k + 1)! (4 k + 3)|
                       | -----                                 |
                       \ k = 0                                 /
\end{maple}

\heading{Conclusion}
\noindent
We have presented an algorithm, and its implementation in \Maple,
to compute FPS. Since it is a goal of Computer Algebra
to work with {\em formal\/} objects, we think this is a very
powerful and valuable addition.

The algorithm is beside of rational operations based on two
basic tools: solving systems of equations and computing
symbolic limits. No truncated series is ever computed.
For the case of asymptotic series, one sided limits are
computed and for all other cases complex ones. 
The power of the procedure for computing LPS stands and falls with
the capabilities of the tool for computing limits.
The algorithms which are used in \Maple\ to compute
limits are described in~\cite{Geddes88,GonnetGruntz92}.

Further, in cases when the resulting RE cannot be solved explicitly,
it can, in principle be used to calculate the coefficients iteratively 
in a lazy evaluation scheme. This is particularly efficient as the
resulting RE always is homogeneous and linear, so that
each coefficient can be calculated by finitely many of its predecessors.
We note that the algorithm to find a simple DE
works, and so such a RE always exists, if
the input function is constructed by integration,
differentiation, addition, multiplication, and the composition with
rational functions or rational powers, from functions with the same 
property, see~\cite{Koediffgl}.  This gives a huge class of functions to 
which the method can be applied.

\end{document}